\documentclass[12pt]{amsart}
\usepackage[utf8]{inputenc}
\usepackage[french]{babel}
\usepackage{geometry} 
\usepackage{graphicx}				
\usepackage{amssymb}
\usepackage{epstopdf}
\usepackage[colorlinks=true,citecolor=cyan,urlcolor=blue]{hyperref} 

\usepackage{abstract}

\newtheoremstyle{mes_theoremes}{1.5em}{2em}{}{}{\bfseries}{~:~}{\parskip}{\thmname{#1}\thmnumber{ #2}\thmnote{ (#3)}}

\theoremstyle{mes_theoremes}

\newtheorem*{pre}{Preuve}

\newtheorem{theo}{Théorème}
\newtheorem{prop}{Proposotion}

\newtheorem*{ex*}{Exemple}
\newtheorem{ex}{Exemple}

\newtheoremstyle{mes_preuves}{1.5em}{2em}{}{}{\it}{~:~}{\parskip}{\thmname{#1}\thmnumber{}}

\theoremstyle{mes_preuves}

\title{Nouvelles conditions pour l'inexistence des nombres parfaits impairs}

\author{\href{mailto:wallace.nancy@courrier.uqam.ca}{Nancy Wallace}}  
\begin{document}
\maketitle

\begin{abstract}

Un nombre, $n$, est dit parfait s'il est \'egal \`a la somme de ses diviseurs propres plus 1. Par exemple $6=1+2+3$. Dans ce document, les deux propositions suivantes seront démontrées: 
\\

S'il existe un nombre parfait impair, $n$, de décomposition en nombre premier $n=p_{1}^{\alpha_{1}}...p_{k}^{\alpha_{k}}q^{\beta}$, où les $ \alpha_{i} $ sont pairs, $ \beta $ est impair et  $q \equiv  5 \mod 8$.  
Alors, au moins un $p_i, ~ 1 \leq i \leq k$ n'est pas un carré dans $\mathbb{Z}/q\mathbb{Z}$. Plus précisément un nombre impair de $p_i$ ne sont pas des carrés dans $\mathbb{Z}/q\mathbb{Z}$.
\\
 
 S'il existe un nombre parfait impair, $n$, de décomposition en nombre premier $n=p_{1}^{\alpha_{1}}...p_{k}^{\alpha_{k}}p_{k+1}^{\beta}$, où les $ \alpha_{i} $ sont pairs, $ \beta $ est impair et $ p_{k+1} \equiv  {1\mod 8}$. 
Alors au moins un $p_i$, $1 \leq i \leq k+1$ est un carré non nul dans au moins un $\mathbb{Z}/{p_j}\mathbb{Z},~ 1 \leq j \leq k+1$.

\end{abstract}

\vspace{25pt}Introduits dans le troisième livre des Éléments de Euclide, les nombres parfaits sont une curiosité mathématique vieille de 5000 ans. Avant de montrer les résultats, voici des exemples permettant de mieux comprendre leurs implications.

\begin{ex}Les nombres $5$, $13$, $29$, $53$ sont tous des carrés dans $\mathbb{Z}/29\mathbb{Z}$. Donc après démonstration de la première proposition, nous serons certains que tout nombre ayant uniquement ces nombres dans sa décomposition en nombres premiers n'est pas un nombre parfait impair si la puissance de 29 est impaire. En principe, cette proposition  est valide seulement avec la contrainte que 5, 13 et 53 ont une puissance paire, mais nous savons par les démonstrations existantes qu'il n'y a qu'un seul nombre dans la décomposition qui est de puissance impaire si $n$ est un nombre parfait impair. Bien entendu, il est possible de construire un ensemble beaucoup plus grand de tels nombres.\end{ex}
\vspace{-30pt}
\begin{ex}Les nombres $5$, $7$, $17$, $73$ n'ont aucun carré non nul dans $\mathbb{Z}/5\mathbb{Z}$, $\mathbb{Z}/7\mathbb{Z}$, $\mathbb{Z}/17\mathbb{Z}$, $\mathbb{Z}/73\mathbb{Z}$. Donc après démonstration da la deuxième proposition, nous serons certains que tout nombre ayant uniquement ces nombres dans sa décomposition en nombres premiers n'est pas un nombre parfait impair et ce peut importe sa puissance. Il est également possible de construire un ensemble contenant une plus grande quantité de tels nombres.
\end{ex}
 Euler a montré qu'un nombre parfait impair, $n$, s'il existe a une décomposition en nombre premier de la forme $n=p_1^{\alpha_1}p_2^{\alpha_2}\cdots p_k^{\alpha_k}q^{\beta}$, où les $\alpha_i, 1\leq i \leq k$ sont paires, $k\geq 2$ et $q \equiv \beta \equiv 1 \mod 4$. Pour les personnes intéressées, une preuve, fait par l'auteure, de ceci est en annexe. Montrons donc le premier résultat annoncé. 

\begin{theo} Si $n$ est un nombre parfait impair de décomposition en nombre premier:
\[n=p_{1}^{\alpha_{1}}...p_{k}^{\alpha_{k}}q^{\beta},\text{ où les $ \alpha_{i} $ sont pairs, $ \beta $ est impair et }  ~q \equiv  5 \mod 8.  \] 
Alors il y a un nombre impair de $p_i, ~ 1 \leq i \leq k$ qui ne sont pas carré dans $\mathbb{Z}/q\mathbb{Z}$.

\begin{pre} Remarquons d'abord que la somme $1+q+...+q^{\beta}$ est paire, car $\beta+1$ est pair, $q$ est impair  et qu'une somme d'un nombre pair de termes impairs est pair. Puisque $n$ est un nombre parfait impair, nous avons :
\[1+q+...+q^{\beta}=2p_{1}^{r_{1}}...p_{k}^{r_{k}}, \text{ où } r_i \leq \alpha_i. \]

Remarquons que la puissance de 2 doit être 1, car autrement $n$ serait pair et nous aurions une contradiction. De plus, $2p_{1}^{r_{1}}...p_{k}^{r_{k}} \equiv 1 \mod q$, alors en posant $q=8m+5$, nous avons, par le symbole de Legendre:
\[\bigg( \frac{2}{q}\bigg)=(-1)^{\frac{q^2-1}{8}}=(-1)^{\frac{(8m+5)^2-1}{8}}=(-1)^{8m^2+10m+3}=((-1)^2)^{4m^2+5m+1}(-1)=-1.\]

Donc:
\begin{align*} \Bigg( \frac{2p_{1}^{r_{1}}...p_{k}^{r_{k}}}{q} \Bigg)&=\bigg( \frac{2}{q}\bigg)\bigg( \frac{p_{1}}{q}\bigg)^{r_{1}}...\bigg(\frac{p_{k}}{q}\bigg)^{r_{k}}=1, \text{ car la formule est multiplicative,}
\\			&=-\bigg( \frac{p_{1}}{q}\bigg)^{r_{1}}...\bigg(\frac{p_{k}}{q}\bigg)^{r_{k}}=1, \text{ car }\bigg( \frac{2}{q}\bigg)=-1.
\end{align*}

Alors en considérant seulement les puissances impaires nous avons:
\[ \bigg( \frac{p_{n_1}}{q}\bigg)^{r_{n_1}}...\bigg(\frac{p_{n_t}}{q}\bigg)^{r_{n_t}}= \bigg( \frac{p_{n_1}}{q}\bigg)...\bigg(\frac{p_{n_t}}{q}\bigg)=-1\text{, où } \{p_{n_j}\} \text{ est une sous suite de } \{p_j\}\]

\vspace{10pt}Nous avons donc bien un nombre impair de puissances impaires qui ne sont pas carrées dans $\mathbb{Z}/q\mathbb{Z}$. En particulier, au moins nombre premier dans la décomposition en nombre premier de $n$ n'est pas un carré dans $\mathbb{Z}/q\mathbb{Z}$ \rule{1.3mm}{1.3mm}
\end{pre}
\end{theo}

Remarquons que cette preuve n'est vraie que si $q\equiv 5 \mod 8$. Nous ferons une preuve similaire avec un résultat différent pour $q\equiv 1 \mod 8$. Avant nous, devons montrer un résultat intermédiaire.

\begin{prop} Si $n$ est un nombre parfait impair de décomposition en nombre premier:
\[n=p_{1}^{\alpha_{1}}...p_{k}^{\alpha_{k}}q^{\beta},\text{ où les $ \alpha_{i} $ sont pairs, $ \beta=\alpha_{k+1} $ est impair et }  ~q= p_{k+1} \equiv  {1\mod 8}.  \] 
Alors il y a un nombre pair de puissances impaires dans la décomposition en nombre premier de $1+p_i+\cdots+p_i^{\alpha_i}$, qui ne sont pas des carrés dans $\mathbb{Z}/{p_i}\mathbb{Z}$. 

\begin{pre} Remarquons d'abord que si $p_i$ est nul dans $\mathbb{Z}/{p_j}\mathbb{Z}$ alors $p_i=p_j$, car $p_i$ et $p_j$ sont premier. 
Puisque $n$ est un nombre parfait, nous avons:
\[ 2n=2p_{1}^{\alpha_{1}}...p_{k}^{\alpha_{k}}p_{k+1}^{\beta}=(1+p_1+...+p_1^{\alpha_1})...(1+p_k+...+p_k^{\alpha_k})(1+p_{k+1}+...+p_{k+1}^{\beta})\]

\vspace{10pt}De plus, pour tout $i$ la somme $1+p_i+...+p_i^{\alpha_i}$ est impaire, car $\alpha_i+1$ et $p_i$ sont impairs et que la somme d'un nombre impair de termes impairs est impair. Puisque la somme divise $2n$ et est impaire, nous avons pour tout $1 \leq i \leq k$:

\[1+p_i+...+p_i^{\alpha_i}=p_{1}^{r_{1}}...p_{i-1}^{r_{i-1}}p_{i+1}^{r_{i+1}}...p_{k}^{r_{k}}p_{k+1}^{r_{k+1}}, \text{ où } r_i \leq \alpha_i, r_{k+1} \leq \beta. \]

\vspace{10pt}De plus, $p_{1}^{r_{1}}...p_{i-1}^{r_{i-1}}p_{i+1}^{r_{i+1}}...p_{k}^{r_{k}}p_{k+1}^{r_{k+1}} \equiv 1 \mod p_i$ et $1$ est toujours un carré dans $\mathbb{Z}/{p_i}\mathbb{Z}$, nous avons alors par le symbole de Legendre:
\[ \Bigg( \frac{p_{1}^{r_{1}}...p_{i-1}^{r_{i-1}}p_{i+1}^{r_{i+1}}...p_{k}^{r_{k}}p_{k+1}^{r_{k+1}}}{p_i} \Bigg)=1.\]

\vspace{10pt}Le symbole de Legendre est une formule multiplicative, nous avons alors:
\[  \bigg( \frac{p_{1}}{p_i}\bigg)^{r_{1}}...\bigg(\frac{p_{i-1}}{p_i}\bigg)^{r_{i-1}}\bigg(\frac{p_{i+1}}{p_i}\bigg)^{r_{i+1}}...\bigg(\frac{p_{k}}{p_i}\bigg)^{r_{k}}\bigg(\frac{p_{k+1}}{p_i}\bigg)^{r_{k+1}}=1.\]

\vspace{10pt}Comme le symbole de Legendre est égal à $1$ ou $-1$, pour tout $r_j$ pair nous avons que $\big( \frac{p_{j}}{p_i}\big)^{r_{j}}=1$ et pour tout $r_j$ impair nous avons que $\big( \frac{p_{j}}{p_i}\big)^{r_{j}-1}=1$. Considérons alors, seulement les $p_j$ ayant une puissance impaire:
\[  \bigg( \frac{p_{n_1}}{p_i}\bigg)...\bigg(\frac{p_{n_t}}{p_i}\bigg)=1 \text{, où } \{p_{n_j}\} \text{ est une sous suite de } \{p_j\}.\]

 D'où il y a un nombre pair de puissances impaires dans la décomposition en nombre premier de $1+p_i+\cdots+p_i^{\alpha_i}$, qui ne sont pas des carrés dans $\mathbb{Z}/{p_i}\mathbb{Z}$.

De façon analogue, la somme $1+p_{k+1}+...+p_{k+1}^{\beta}$ est paire, car $\beta+1$ est pair, $p_{k+1}$ est impair  et qu'une somme d'un nombre pair de termes impair est pair. Et puisque la somme divise $2n$, car $n$ est parfait, nous avons:
\[1+p_{k+1}+...+p_{k+1}^{\beta}=2p_{1}^{r_{1}}...p_{k}^{r_{k}}, \text{ où } r_i \leq \alpha_i. \]
Remarquons que la puissance de 2 doit être 1, car autrement $n$ serait pair et nous aurions une contradiction.
\\ 

Puisque $p_{k+1} \equiv 1 \mod 8$, par supposition de départ, posons  $p_{k+1}=8m+1$. Nous avons alors:
\[ \bigg( \frac{2}{p_{k+1}}\bigg)=(-1)^{\frac{p_{k+1}^2-1}{8}}=(-1)^{\frac{(8m+1)^2-1}{8}}=(-1)^{8m^2+2m}=((-1)^2)^{4m^2+1m}=1. \]

Alors le même argument s'applique également à $p_{k+1}$ et nous avons le résultat voulu \rule{1.3mm}{1.3mm}
\end{pre}
\end{prop}

\begin{theo} Si $n$ est un nombre parfait impair de décomposition en nombre premier:
\[n=p_{1}^{\alpha_{1}}...p_{k}^{\alpha_{k}}p_{k+1}^{\beta},\text{ où les $ \alpha_{i} $ sont pairs, $ \beta $ est impair et }  ~ p_{k+1} \equiv  {1\mod 8}.  \] 
Alors au moins un $p_i$ est un carré non nul dans au moins un $\mathbb{Z}/{p_j}\mathbb{Z},~ 1 \leq j \leq k+1$.

\begin{pre}
Par contradiction, supposons qu'il n'y ait aucun des $p_i$ qui soient un carré non nul dans au moins un $\mathbb{Z}/{p_j}\mathbb{Z},~ 1 \leq j \leq k+1$. 
\\

Pour $1+p_i+...+p_i^{\alpha_i}=p_{1}^{r_{1,i}}...p_{i-1}^{r_{i-1,i}}p_{i+1}^{r_{i+1,i}}...p_{k}^{r_{k,i}}p_{k+1}^{r_{k+1,i}}$, posons $\{ r_{d_1},..., r_{d_s} \}$ l'ensemble des puissances paires et $\{ r_{n_1},..., r_{n_t} \}$ l'ensemble  des puissances impaires. Nous avons alors que $t$ est un indice pair par la proposition précédente, car il n'y a aucun carré par hypothèse.
Donc la somme $r_{1,i}+...+r_{k+1,i}= r_{d_1}+...+r_{d_s}+ r_{n_1}+...+r_{n_t}$ est paire, car elle contient un nombre pair de nombres impairs. Nous trouvons le même résultat pour $1+p_{k+1}+...+p_{k+1}^{\beta}$.
\\

Nous avons alors pour:
\begin{align*}2n&=(1+p_1+...+p_1^{\alpha_1})...(1+p_k+...+p_k^{\alpha_k})(1+p_{k+1}+...+p_{k+1}^{\beta})
\\			&=(p_{2}^{r_{2,1}}...p_{k+1}^{r_{k+1,1}})...(p_{1}^{r_{1,k}}...p_{k-1}^{r_{k-1,k}}p_{k+1}^{r_{k+1,k}})(2p_{1}^{r_{1,k+1}}...p_{k}^{r_{k,k+1}})
\end{align*}
\[ \text{ que la somme } (r_{2,1}+...+r_{k+1,1})+...+ ( r_{1,k}+...+r_{k-1,k}+r_{k+1,k})+(r_{1,k+1}+...+r_{k,k+1})\] 
est paire, car chaque parenthèse est paire par ce qui précède.
\\

Or $ r_{1,2}+...+r_{1,k+1}=\alpha_1,~ ..., r_{k,1}+...+r_{k,k-1}+r_{k,k+1}=\alpha_k,~ r_{k+1,1}+...+r_{k+1,k}=\beta$ et nous avons que $\alpha_1+...+\alpha_k+\beta $ est  impaire, car les $\alpha_i$ sont tous pairs et $\beta$ est impair. D'où la contradiction \rule{1.3mm}{1.3mm}
\end{pre}

\end{theo}

Remarquons que le théorème 2 est moins restrictif que le premier. Par exemple, si deux nombres premiers, disons $p_1$ et $p_2$ sont tels que $p_1\equiv p_2 \equiv 3 \mod 4$ et que $p_1$ n'est pas un carré dans $p_2$, alors $p_2$ est un carré dans $p_1$. Donc il suffit que deux parmi les $\{p_i|1\leq i\leq k\}$ soient congru à 3 modulo 4 et la condition est satisfaite.

En effet, soit $p_1=4n+3,~ p_2 =4m+3$ tel que $p_1$ ne soit pas un carré dans $p_2$. Alors par la formule de réciprocité quadratique de Gauss-Legendre, nous avons:
\begin{align*} \bigg( \frac{p_2}{p_1} \bigg) &=(-1)^{\frac{(4n+3-1)(4m+3-1)}{4}}\bigg( \frac{p_1}{p_2} \bigg)=(-1)^{\frac{16mn+8m+4n+4}{4}}\bigg( \frac{p_1}{p_2} \bigg)
\\								 &=(-1)^{4mn+2m+2n+1}(-1)=(-1)^{4mn+2m+2n+2}=((-1)^2)^{2mn+m+n+1}=1
\end{align*}
Donc $p_2$ est bien un carré dans $p_1$.

Comme Jacques Lefèvre, l'auteure croit qu'il n'existe pas de nombres parfaits impairs et ces deux résultats sont un bon pas dans la bonne direction afin de montrer leurs inexistences. De plus, pour les personnes faisant un programme cherchant un nombre parfait impair, ceci donne des critères supplémentaires pour éviter de calculer inutilement certains nombres.

\newpage

\begin{center}\begin{Huge} Annexe \end{Huge}\end{center}

\begin{align*}\text{D\'efinition :} :& \text{ Un nombre $n$ est dit $parfait$ s'il est \'egal \`a la somme de ses diviseurs propres. }
\\						& \text{ Dans cet annexe 1 sera consid\'er\'e comme un diviseur propre.}
\end{align*}

Notation : Nous noterons $n'=\underset{d|n,d<n}{\sum}d$

\begin{align*} \text{Th\'eor\`eme :} &\text{ Si $n$ est un nombre parfait impair, alors sa d\'ecomposition en nombre premier}
\\						     &\text{ satisfait : } n=p_{1}^{\alpha_{1}}...p_{k}^{\alpha_{k}}q^{\beta},
\\						     & \text{ o\`u les } \alpha_{i} \text{ sont pairs,  } k \geq 2, ~ \beta \text{ impair et } ~ q \equiv \beta \equiv 1 \mod~4.
\end{align*}

\begin{align*} \text{Lemme 0 :}  &\text{ Si $n$ est un nombre parfait impair, alors } n \neq p^\alpha \text{,  $p$ premier.}
\\			
\\			\text{Preuve : }  &\text{ Remarquons d'abord que les diviseurs propres de $n$ sont : } 1,p,...,p^{\alpha -1 }
\\						&\hspace{80pt}n=p^\alpha=(p-1)(p^{\alpha -1 }+...+p+1)+1,
\\						&\hspace{100pt}\Rightarrow n > p^{\alpha -1 }+...+p+1  ~\rule{1.3mm}{1.3mm}
\end{align*}

\begin{align*} \text{Lemme 1 :}  &\text{ Si } n=p_{1}^{\alpha_{1}}...p_{k}^{\alpha_{k}} \text{ o\`u les $p_i$ sont premiers, alors $n$ a $(\alpha_{1}+1)...(\alpha_{k}+1)$ diviseurs. }
\\						
\\			\text{Preuve :} &\text{ Tout diviseur de $n$ est un produit de la forme } p_{1}^{\beta_{1}}...p_{k}^{\beta_{k}}, \text{ o\`u } 0 \leq \beta_i \leq \alpha_i ~\forall ~i.
\\						&\text{ il existe donc }(\alpha_{i}+1) \text{ possibilit\'es pour chaque } ~\beta_i.
\\  
\\						&\text{ Donc d'après l'unicité de la factorisation en nombres premiers,}
\\						&\text{ $n$ a } ~(\alpha_{1}+1)...(\alpha_{k}+1)\text{ diviseurs distincts } ~\rule{1.3mm}{1.3mm}
\end{align*}

\begin{align*} \text{Lemme 2 :}  &\text{ Si $n$ est un nombre parfait impair, alors sa d\'ecomposition en nombres premiers}
\\						 &\text{ satisfait : } n=p_{1}^{\alpha_{1}}...p_{k}^{\alpha_{k}}, \text{ o\`u  $k \geq 2$ et les } \alpha_{i} \text{ ne sont pas tous pairs.  }
\\
\\			\text{Preuve : }  &\text{ Par le lemme 1, $n$ a }~ (\alpha_{1}+1)...(\alpha_{k}+1)-1\text{ diviseurs propres, car $n$ est un}
\\						&\text{ diviseur impropre de $n$. Ainsi en supposant que tous les $\alpha_{i}$ dans la d\'ecomposition} 
\\						&\text{ en nombres premiers de $n$ sont pairs, nous aurons que $n'$ est une somme d'un  }
\\						&\text{ nombre pair de termes impairs et donc que $n'$ est pair. }
\\
\\						&\text{ Donc $n$ ne peut pas \^etre un nombre parfait puisqu'il est impair  }  \rule{1.3mm}{1.3mm}
\end{align*}

\begin{align*} \text{Lemme 3 :}  &\text{ Si $n$ est un nombre parfait impair, alors sa d\'ecomposition en nombres premiers}
\\						 &\text{ satisfait : } n=p_{1}^{\alpha_{1}}...p_{k}^{\alpha_{k}}q^{\beta}, \text{ o\`u  les } \alpha_{i} \text{ sont pairs et } \beta \text{ est impair.}
\\
\\			\text{Preuve : }  	&\text{ Par contradiction, supposons que $n$ est un nombre parfait et que sa d\'ecomposition }
\\						&\text{ en nombres premiers satisfasse : }n=p_{1}^{\alpha_{1}}...p_{k}^{\alpha_{k}}q_{1}^{\beta_{1}}...q_{m}^{\beta_{m}},\text{ o\`u  les } \alpha_{i} \text{ sont pairs et les } \beta_i  
\\						&\text{ sont impairs et avec m $>$1. }  \text{ Nous aurions alors:}
\\						&\hspace{60pt}2n=(\sum_{i=0}^{\alpha_{1}} p_{1}^i)...(\sum_{i=0}^{\alpha_{k}} p_{k}^i)\times \underbrace{(\sum_{i=0}^{\beta_{1}} q_{1}^i)}_{\beta_1 + 1 ~ \text{termes}}...\underbrace{(\sum_{i=0}^{\beta_{m}} q_{m}^i)}_{\beta_1 + 1 ~ \text{termes}},
\\						&\text{ Donc la somme de chacune des parenth\`eses à droite du signe multiplicatif $\times$ est  }
\\						&\text{ paire, car elle contient une somme d'un nombre pair de termes impairs. Posons :}
\\						&\hspace{90pt} 1+q_{i}+...+q_{i}^{\beta_{i}}=2z_{i}, ~ z_{i} \in \mathbb{N}
\\						&\hspace{100pt}\Rightarrow 2n=(\sum_{i=0}^{\alpha_{1}} p_{1}^i)...(\sum_{i=0}^{\alpha_{k}} p_{k}^i)2z_{1}...2z_{m},
\\						&\hspace{130pt}= (\sum_{i=0}^{\alpha_{1}} p_{1}^i)...(\sum_{i=0}^{\alpha_{k}} p_{k}^i)2^mz_{1}...z_{m}
\\
\\						&\text{ Ce qui est absurde, car $m > 1$ et $n$ est impair. }
\\						& \text{ Donc $m \leq 1$ et par le lemme 2, $m=1$ }\rule{1.3mm}{1.3mm}
\end{align*}
\\

\begin{align*} \text{Lemme 4 :} &\text{ Si $n$ est un nombre parfait impair, alors sa d\'ecomposition en nombres }
\\						&\text{ premiers satisfait : }n=p_{1}^{\alpha_{1}}...p_{k}^{\alpha_{k}}q^{\beta},\text{ où les $\alpha_{i}$ sont pairs et } \beta \equiv 1 ~mod ~4.
\\			
\\			\text{Preuve :}   &\text{ Supposons que $n$ est un nombre parfait de d\'ecomposition en nombre premier}
\\						&n=p_{1}^{\alpha_{1}}...p_{k}^{\alpha_{k}}q^{\beta},\text{ où les $\alpha_{i}$ sont pairs et } \beta \equiv 3~ mod~ 4.
\\						&\text{ Nous avons, en posant $\beta+1=4m$ \big(car $4 \Big| \beta+1 \big)$ et $q=2^{x}z+1$ avec $x \geq 1, z$ impair :}
\end{align*}
\begin{align*}						 \Bigg(\sum_{i=0}^{\beta} &q^i \Bigg)=\Bigg(\frac{q^{\beta+1}-1}{q-1} \Bigg)= \Bigg(\frac{({2^{x}z+1})^{4m}-1}{2^{x}z+1-1} \Bigg)
\\
\\						= &\frac{\binom{4m}{0} 2^{4mx}z^{4m}+\binom{4m}{1} 2^{(4m-1)x}z^{4m-1}+...+\binom{4m}{4m-3} 2^{3x}z^3+\binom{4m}{4m-2} 2^{2x}z^2+\binom{4m}{4m-1} 2^{x}z+1-1}{2^xz}
\\
\\						=&2^{(4m-1)x}z^{4m-1}+4m2^{(4m-2)x}z^{4m-2}+...+\binom{4m}{4m-3}2^{2x}z^2+\frac{4m(4m-1)}{2} 2^xz+4m
\\
\\						=&2^{(4m-1)x}z^{4m-1}+4m2^{(4m-2)x}z^{4m-2}+...+\binom{4m}{4m-3}2^{2x}z^2+2m(4m-1) 2^xz+4m
\\
\\						= &4(2^{(4m-1)x-2}z^{4m-1}+4m2^{(4m-2)x-2}z^{4m-2}+...+\binom{4m}{4m-3}2^{2x-2}z^2+m(4m-1) 2^{x-1}z+m ), 
\\						&\text{ car $x \geq 1$ et $m \geq 1$. D'où :}
\\						& \hspace{150pt} \Rightarrow 4 \Big| \sum_{i=0}^{\beta} q^i 
\\						&\hspace{90pt} \Rightarrow 4 \Big| \Big(\sum_{i=0}^{\alpha_{1}} p_{1}^i\Big)...\Big(\sum_{i=0}^{\alpha_{k}} p_{k}^i\Big)\Big(\sum_{i=0}^{\beta} q^i\Big)=2n
\\
\\						& \text{ Ce qui est absurde, car $n$ est impair }\rule{1.3mm}{1.3mm}
\end{align*}
\newpage
\begin{align*} \text{Lemme 5 :} &\text{ Si $n$ est un nombre parfait impair, alors sa d\'ecomposition en nombres }
\\						&\text{ premiers satisfait : }n=p_{1}^{\alpha_{1}}...p_{k}^{\alpha_{k}}q^{\beta},\text{ où les $\alpha_{i}$ sont pairs, } \beta \equiv 1 ~mod ~4 
\\						&\text{  et  } q \equiv 1 ~mod ~4 .
\\			
\\			\text{Preuve :}   &\text{ Supposons que $n$ est un nombre parfait de d\'ecomposition en nombres premiers}
\\						&n=p_{1}^{\alpha_{1}}...p_{k}^{\alpha_{k}}q^{\beta},\text{ tel que les $\alpha_{i}$ sont pairs, } \beta \equiv 1~ mod~ 4 \text{  et  } q \equiv 3 ~mod ~4 .
\\						&\text{ En posant $\beta=4m+1$ et $q=4z-1$, nous avons :}
\\
\\						& q^{\beta+1}=(4z-1)^{4m+2}=(4z-1)^2(4z-1)^{4m},
\\
\\						&\hspace{22pt}=16z^2((4z-1)^{4m}-8z(4z-1)^{4m}+(4z-1)^{4m},
\\
\\						& \Rightarrow \sum_{i=0}^{\beta} q^i =\frac{q^{\beta+1}-1}{q-1} =\frac{16z^2((4z-1)^{4m}-8z(4z-1)^{4m}+(4z-1)^{4m}-1}{4z-2},
\\
\\						&\hspace{46pt}=\frac{4z(4z-1)^{4m}(4z-2)+(4z-1)^{4m}-1}{4z-2},
\\
\\						&\hspace{46pt}=4z(4z-1)^{4m}+\frac{(4z-1)^{4m}-1}{4z-2},
\\
\\						&\hspace{46pt}=4z(4z-1)^{4m}+\frac{\binom{4m}{0} 4^{4m}z^{4m}+...+\binom{4m}{4m-2} 4^2z^2+\binom{4m}{4m-1} 4z+1-1}{4z-2},
\\
\\						&\hspace{46pt}=4z(4z-1)^{4m}+\frac{\binom{4m}{0} 4^{4m}z^{4m}+...+\binom{4m}{4m-2} 4^2z^2+4m 4z+1-1}{4z-2},
\\
\\						&\hspace{46pt}=4z(4z-1)^{4m}+\frac{16( 4^{4m-2}z^{4m}+...+\binom{4m}{4m-2} z^2+mz)}{4z-2},
\\
\\						&\hspace{46pt}=4z(4z-1)^{4m}+\frac{8( 4^{4m-2}z^{4m}+...+\binom{4m}{4m-2} z^2+mz)}{2z-1}.
\\
\end{align*}
	 Puisque nous savons que $\frac{q^{\beta+1}-1}{q-1}  \in \mathbb{N}$ et que $4z(4z-1)^{4m} \in \mathbb{N}$, nous en déduisons que :						
	\[ \frac{8(\binom{4m}{0} 4^{4m-2}z^{4m}+...+mz)}{2z-1} \in \mathbb{N}.\]
	De plus, $2z-1\nmid 2$ donc :
	\[ \frac{4^{4m-2}z^{4m}+...+mz}{2z-1} \in \mathbb{N}.\]
	
	\[\text{Comme } \frac{q^{\beta+1}-1}{q-1} =4\Bigg(z(4z-1)^{4m}+\frac{2( 4^{4m-2}z^{4m}+...+mz)}{2z-1}\Bigg),\]
	
	\[\text{nous avons que } 4 \Bigg| (\sum_{i=0}^{\alpha_{1}} p_{1}^i)...(\sum_{i=0}^{\alpha_{k}} p_{k}^i)(\sum_{i=0}^{\beta} q^i)=2n \]
\\
\\	
	Ce qui est absurde, car $n$ est impair. \rule{1.3mm}{1.3mm}
\\
\\


\begin{align*} \text{Lemme 6 :} &\text{ Si $n$ est un nombre parfait impair, alors sa d\'ecomposition en nombres }
\\						&\text{ premiers ne satisfait pas: }n=p^{\alpha}q^{\beta},\text{ où $\alpha$ est pairs et $\beta$ est impair. } 
\\
\\			\text{Preuve :}  &\text{ Supposons que $n$ soit un nombre parfait impair, et que sa d\'ecomposition  }
\\						&\text{ en nombre premier satisfasse: }n=p^{\alpha}q^{\beta},\text{ où $\alpha$ est pair, $\beta$ est impair et $\alpha, \beta \geq 1$. } 
\\	
\\						&\text{ Remarquons d'abord que } p,q\geq 3 \text{, car $n$ est impair.}
\\						& \text{ De plus, nous avons:}
\\
\\						&\hspace{100pt} p \nmid  p^{\alpha}+p^{\alpha-1}+...+1  ~\text{ et  }~  q \nmid q^{\beta} +q^{\beta-1}+...+1.
\\
\\						&\text{ Nous avons aussi que } p^{\alpha}+p^{\alpha-1}+...+1  \text{ est impair, car c'est une somme d'un nombre }
\\						&\text{ impair de nombres impairs.}
\\
\end{align*}
\begin{align*}
\\						&\text{ Nous avons donc : }
\\
\\						& \hspace{100pt}2n= \Bigg(\frac{p^{\alpha+1}-1}{p-1} \Bigg) \Bigg(\frac{q^{\beta+1}-1}{q-1} \Bigg)
\\
\\						&\hspace{120pt} \Rightarrow \frac{q^{\beta+1}-1}{q-1}  \text{ est pair.}
\\
\\					    (1)& \hspace{40pt}\Rightarrow \frac{p^{\alpha+1}-1}{p-1} =q^{\beta} \text{ et } \frac{q^{\beta+1}-1}{q-1} =2p^{\alpha} \text{, car }n=p^{\alpha}q^{\beta},
\\
\\				    	   (2) & \hspace{70pt} \Rightarrow   \frac{p^{\alpha+1}-1}{p-1} =q^{ \beta } \Rightarrow p^{\alpha+1}=q^{\beta}p-q^{\beta}+1,
\\
\\				   	    (3)&\hspace{10pt}  \text{ et }   \frac{q^{\beta+1}-1}{q-1} =2p^{ \alpha } \Rightarrow q^{\beta+1}=2p^{\alpha}(q-1)+1=2p^{\alpha}q-2p^{\alpha}+1,
\\
\\			                            &  \hspace{50pt}\text{ et }2p^{\alpha}-1=\frac{q^{\beta+1}-1}{q-1} -1=\frac{q^{\beta+1}-q}{q-1} =q\frac{q^{\beta}-1}{q-1},
\\
\\						&\hspace{30pt} \text{ et }  2p^{\alpha}-1=p^{\alpha}+p^{\alpha}-1=p^{\alpha}+(p-1)(p^{\alpha-1}+...+p+1).
\\
\\
\\						&\text{ Les deux dernières implications nous donnent :}
\\						& \hspace{46pt}q\frac{q^{\beta}-1}{q-1}=p^{\alpha}+(p-1)(p^{\alpha-1}+...+p+1),
\\
\\						&\hspace{56pt}\Rightarrow q \mid p^{\alpha}+(p-1)(p^{\alpha-1}+...+p+1),
\\
\\			        	             (4)&\hspace{70pt}\text{ et } q \nmid p^{\alpha-1}+...+p+1 \text{, car } q \nmid p^{\alpha}.
\\
\end{align*}
\begin{align*}
\\						& \text{ Puisque } q \mid p^{\alpha}+(p-1)(p^{\alpha-1}+...+p+1) \text{ et que par (1) } q \mid p^{\alpha}+p^{\alpha-1}+...+p+1, 
\\						&\text{ nous avons donc :}
\\
\\						& \hspace{46pt}q \mid p^{\alpha}+(p-1)(p^{\alpha-1}+...+p+1)-(p^{\alpha}+p^{\alpha-1}+...+p+1),
\\						& \hspace{66pt}q \mid (p-1)(p^{\alpha-1}+...+p+1)-(p^{\alpha-1}+...+p+1),
\\						&\hspace{106pt}\Rightarrow q \mid (p-2)(p^{\alpha-1}+...+p+1),
\\						&\hspace{50pt}\Rightarrow q \mid p-2, \text{ car par (4) } ~ q \nmid p^{\alpha-1}+...+p+1 , 
\\
\\			       \text{(5)  }&\text{ Il existe donc $s \in \mathbb{N}$ tel que } sq=p-2.\text{ Notons que } s \geq 1 \text{, car } p \geq 3.
\\
\\
\\						&\text{ De plus, nous avons :}
\\
\\						& \hspace{70pt}p^{\alpha+1}q=q^{\beta+1}p-q^{\beta+1}+q, \text{ par (2)},
\\						& \hspace{98pt}=(2p^{\alpha}(q-1)+1)p-2p^{\alpha}(q-1)-1+q , \text{ par (3)},
\\						& \hspace{98pt}=2p^{\alpha+1}q-2p^{\alpha+1}+p-2p^{\alpha}q+2p^{\alpha}-1+q , 
\\						& \hspace{98pt}=p(2p^{\alpha}q-2p^{\alpha}+1-2p^{\alpha-1}q+2p^{\alpha-1})-1+q , ,\text{ car } \alpha \geq 1
\\
\\						&\hspace{70pt} \Rightarrow p \big| p(2p^{\alpha}q-2p^{\alpha}+1-2p^{\alpha-1}q+2p^{\alpha-1})-1+q , 
\\						& \hspace{125pt}\Rightarrow p \big|q-1. 
\\
\\						&\text{ Il existe donc $t \in \mathbb{N}$ tel que } tp=q-1.\text{ Notons que } t \geq 1 \text{, car } q \geq 3.
\\						& \text{ Nous avons alors, par (5):}
\\
\\						& \hspace{106pt} sq+2=p \Rightarrow sqt+2t=q-1
\\						&\hspace{126pt}\Rightarrow sqt-q=-2t-1. 
\\
\\						&\text{ Ce qui est absurde, puisque $s$ et $t$ sont positifs et donc } sqt-q \geq 0 \text{ et } -2t-1 \leq -1.
\\
\\						&\text{ D'o\`u si $n$ est un nombre parfait impair, alors sa d\'ecomposition en nombres }
\\						&\text{ premiers ne satisfait pas: }n=p^{\alpha}q^{\beta},\text{ où $\alpha$ est pair et $\beta$ est impair} ~\rule{1,3mm}{1.3mm} 
\\
\\	\text{ Preuve }			 &\text{ du th\'eor\`eme : La preuve du th\'eor\`eme d\'ecoule directement des lemmes 0,2,3,4,5,6.}~\rule{1,3mm}{1.3mm}
\end{align*}


\begin{thebibliography}{99}
\bibitem{1}  \begin{sl}\'Eléments\end{sl},  Euclide, livre IX, III$^e$ si\`ecle avant J.C.
\bibitem{2}  {\begin{sl} Le Tractatus de numerorum doctrina capita sedecim, quae supersunt. \end{sl}}, Leonhard Euler, Commentationes arithmeticae 2, (1849), pp. 503-575
\bibitem{3}   {\begin{sl}Th\'eorie alg\'ebrique des nombres.\end{sl}},P. Samuel,  \'Edition Hermann, collection Méthodes, (1967).
\bibitem{4} {\begin{sl}History of the theory of numbers. Vol. I: Divisibility and primality.\end{sl}}, L. E.  Dickson,  Publi\'e par Carnegie Institution of Washington, (1919),  520 pages.
\end{thebibliography}
\end{document}